\documentclass[12pt,leqno]{article}
\usepackage[francais]{babel} 
\usepackage[latin1]{inputenc}
\usepackage[T1]{fontenc}
\usepackage{amsmath}
\usepackage{amsfonts}
\usepackage{amssymb}

\newcommand{\C}{\mathbb{C}}
\newcommand{\R}{\mathbb{R}}

\newcommand{\Z}{\mathbb{Z}}
\newcommand{\1}{1\!\!\!1}

\newcommand{\Hh}{\mathcal{H}}

\newcommand{\intint}{\int\!\!\!\!\int}

\begin{document}
\title{Sur les séries de Fourier des fonctions continues unimodulaires}
\author{Jean Bourgain et Jean--Pierre Kahane}
\maketitle

Haim Brézis a découvert une formule intéressante qui donne le degré topologique d'une fonction continue unimodulaire en fonction de ses coefficients de Fourier ; à savoir
\begin{equation}
\deg f = \sum_{-\infty}^{\infty} n|a_n|^2
\end{equation}
sous les hypothèses ($ \sim $ se lit \og a pour série de Fourier\fg)
\begin{equation}
f(e^{it})\sim \sum_{-\infty}^{\infty} a_n e^{int}, \quad f\in C(S^1,S^1)\,,
\end{equation}
\begin{equation}
\hskip 4cm\sum_{-\infty}^\infty |n|\, |a_n|^2 < \infty\,. \hskip5cm [2]
\end{equation}
Rappelons que $\deg f=k$ quand $f(z)=z^k$ $(k\in \Z)$ et que $\deg f=\deg g$ quand $f$ et $g$ sont homotopes dans $C(S^1,S^1)$.
L'hypoth\`ese (3) signifie $ f\in H^{1/2}(S^{1},S^{1}) $. A partir de l\`{a} Br\'{e}zis a pos\'{e} la question : est--il vrai que (2) entra\^{i}ne toujours
\begin{equation}
\sum_{-\infty}^{\infty}|n|\, |a_n|^2 \le |\deg f| + \sum_{0}^{\infty} n |a_n|^2\ \ ?
\end{equation} 
De fa\c{c}on \'{e}quivalente, est--il vrai que sous l'hypoth\`{e}se (2) on ait l'implication
\begin{equation}
\hskip2cm\sum_0^\infty n |a_n|^2 < \infty \Longrightarrow  \sum_{-\infty}^{\infty} |n|\, |a_n|^2 <\infty\ ? \hskip3cm [3]
\end{equation}

Cette simple question montre bien que l'analyse harmonique des fonctions unimodulaires est un sujet riche et peut r\'{e}v\'{e}ler des ph\'{e}nom\`{e}nes nouveaux et int\'{e}ressants.

Nous allons \'{e}largir la question et y r\'{e}pondre positivement.

\vskip2mm

\textsc{Th\'{e}or\`{e}me 1}.-- \textit{Soit $ 0<s<1$. Sous l'hypoth\`{e}se $ (2) $ on a l'implication}
\begin{equation}
\sum_0^\infty n^{2s}|a_n|^2 <\infty \Longrightarrow \sum_{-\infty}^\infty |n|^{2s} |a_n|^2 <\infty\,.
\end{equation}

Le membre de droite signifie $ f\in H^s(S^1,S^1)$. Quitte \`{a} décaler les coefficients, nous supposerons $ \deg f=0 $ et nous montrerons que l'hypoth\`{e}se 
\begin{equation}
\sum_0^\infty n^{2s}|a_n|^2 \le C <\infty
\end{equation}
entra\^{i}ne
\begin{equation}
\sum_{-\infty}^\infty |n|^{2s}|a_n|^2 \le C'(f) <\infty\,.
\end{equation}

Br\'{e}zis et Nirenberg ont \'{e}tendu la notion de degr\'{e} topologique aux fonctions de la classe $VMO$, c'est--\`{a}--dire limites de fonctions continues dans la norme de $BMO$ (voir [4]). Rappelons la d\'{e}finition :
\begin{equation}
\parallel g\parallel_{BMO(S^1,S^1)}= \sup_{t\in \R,\, s>0} \frac 1{2s} \int^s_{-s} \Big| g(e^{i(t+u)})- \frac 1{2s} \int^s_{-s} g(e^{i(t+u')} du'\Big|
du.
\end{equation}
On peut ainsi étendre le théorème 1.

\vskip2mm

\textsc{Théoreme 2}.-- \textit{Soit toujours $0<s<1$. L'implication $(6)$ est valable sous l'hypothèse}
\begin{equation}
f(e^{it}) \sim \sum_{-\infty}^\infty a_n e^{int},\ f\in VMO(S^1,S^1)\,.
\end{equation}

Avant de passer aux démonstrations, examinons deux questions :

\vskip2mm

\textbf{Q1.} Peut--on, dans les hypothèses du théorème 1, remplacer $f\in C(S^1,S^1)$ par $f\in L^\infty(S^1,S^1)$ ?

\vskip 2mm

\textbf{Q2.} Peut--on, pour obtenir l'implication $(7)\Longrightarrow (8)$, remplacer $f\in C(S^1,S^1)$ 
ou $f\in VMO(S^1, S^1)$ par une hypoth\`ese du type $f\in H^{s'} (S^1,S^1)$ ?

Les réponses sont négatives (en ce qui concerne Q2, pour $s' <\frac{1}{2}$)  et fondées sur les produits de Blaschke.

\vskip2mm

\textbf{R1}.
 Soit $g(z)$ un produit de Blaschke tel que $g(0)=0$. Posons $f(e^{it}) = g(e^{-it})$. Alors $f\in L^\infty (S^1,S^1)$ et le  premier membre de (7) est nul. 
 Montrons maintenant qu'on peut choisir $g$ de sorte que $g(e^{it})$ n'appartienne à aucun espace $H^s$ $(s>0)$ ; cela achèvera la preuve que la réponse à la question Q1 est négative. Plus généralement, montrons qu'étant donné une suite positive $\omega=(\omega_n)_{n\in \mathbb{N}}$ tendant vers l'infini, il existe un produit de Blaschke $B(z)=\sum\limits_0^\infty b_nz^n$ tel que $\sum\limits_0^\infty |b_n|^2\omega_n=\infty$. En effet, partons d'un produit de Blaschke infini que nous écrivons $\prod\limits_1^\infty B_j(z)$, où les $B_j(z)$ sont des produits de Blaschke finis. Il existe alors une suite d'entiers $\nu_j$ tendant vers l'infini telle que $B(z) =\prod\limits_1^\infty B_j(z^{\nu_j})$ ait la propriété voulue. Pour la construire, imposons la condition que $\nu_j$ divise $\nu_{j+1}$ $(i=1,2,\ldots)$. Les sommes partielles d'ordre $ \nu_{k+1} $ des séries de Taylor de $ B(z) $ et de $ C_k(z)=\prod\limits_1^k B_j (z^{\nu_j}) $ sont les mêmes à un facteur multiplicatif près qui tend vers 1 quand $ k\rightarrow\infty $. Appelons norme d'un produit de Blaschke et notons $\|\ \|$ la norme de la suite de ses coefficients de Taylor dans $\ell^2(\mathbb{N},\omega)$. Pour avoir la propriété voulue, à savoir $\|B\|=\infty$, il suffit que les normes des sommes partielles d'ordre $\nu_{k+1}$ des $C_k$ tendent vers l'infini $(k\rightarrow\infty)$ ; il suffit que $\|C_k\|\rightarrow\infty$ (réalisable si les $\nu_k$ croissent assez vite) et que les $\nu_{k+1}$ croissent assez vite (conditions réalisables par induction).
 
 \vskip2mm
 
 \textbf{R2.} Revenons au produit de Blaschke $g(z)$, et prenons $f(e^{it})=e^{-it}\break g(e^{-it})$. Si $f\in VMO(S^1,S^1)$, et en particulier si $f\in H^{1/2}(S^1,S^1)$, le degré topologique de $g$ est fini, donc $g(z)$ est une fraction rationnelle. Prenons maintenant pour $g(z)$ un produit de Blaschke qui n'est pas une  fraction rationnelle, et dont les coefficients de Taylor sont $O\Big(\frac{1}{n}\Big)$ ($n\rightarrow \infty)$ (Newman et Shapiro 1962 [4]) ; 
alors $f\in H^{s'}(S^1,S^1)$ pour tout $s'<\frac{1}{2}$, on peut prendre $C=0$ dans (7) et le premier membre de (8) est infini pour $s\geq \frac 12$. 
La réponse à Q2 est donc négative quand $s'<\frac{1}{2}$. Quand $s'>\frac{1}{2}$ on est ramené à $f\in C(S^1,S^1)$, et pour $s'=\frac{1}{2}$ au cas étudié par Brézis, donc la remarque est positive.

\vskip4mm
Le referee s'est demandé où intervient l'hypothèse $s<1$ dans la démonstration du théorème~1. Nous répondons à cette question dans la phrase suivant (46). Mais cela suggère une nouvelle question :

\vskip2mm

\textbf{Q3.}= Peut--on étendre les théorèmes 1 et 2 en remplaçant l'hypothèse $0<s<1$ par $s>0$ ?

Nous verrons à la fin de l'article que la réponse est positive : c'est le théorème~3.
 
La preuve du théorème 1 se fera en trois temps : $ s=\frac{1}{2}$ (solution de la question de Brézis), $ s>\frac{1}{2} $ et $ s<\frac{1}{2} $. La preuve du théorème 2 en sera une adaptation.

\vskip2mm
\textbf{Remarque.}

On peut se poser la question si l'in\'egalit\'e (8) reste valable avec $C'(C)$ (en supposant deg\,$f=0$) comme c'est le cas pour $s=\frac 12$.
La r\'eponse est n\'egative comme le montre l'exemple suivant. Posons, pour $0<a<1$ et $k\in\mathbb Z_+$
\begin{equation*}
f(e^{it})= e^{-ikt} \frac{a-e^{ikt}}{1-ae^{ikt}}=\frac 1a e^{-ikt}-(1-a^2)\sum_{j\geq -1} a^j e^{ijkt}
\end{equation*}

Alors deg\,$f=0$.

Laissons $a\to 1$.
Pour $0<s<\frac 12$, la norme $H^s$ du second terme est de l'ordre de $(1-a)^{\frac 12-s}k^s\to 0$ pour $k$ fix\'e (donc $C$ est arbitrairement petite)
tandis que $(8)\sim k^{2s}$.
Pour $s>\frac 12$, posons $k=1$ et consid\'erons la fonction $\overline{f(e^{it})}$.
Alors $C\sim 1$ dans (7) et $C'\sim\Big(\frac 1{1-a}\Big)^{2s-1}$ dans (8).
Par contre, en supposant que $\Vert f\Vert_{BMO}$ est suffisamment petite, on a bien $C'=C'(C)$, comme on le montrera dans la preuve du
Th\'eor\`eme 2.

\vskip2mm

\textbf{Démonstration du théorème }1.

\textbf{Cas } $ s=\frac{1}{2}  $. Nous nous proposons d'établir (8) à partir de (2), (7) et $ \deg f=0 $.

Soit $ K_\varepsilon(\varepsilon\downarrow 0) $ une suite de noyaux régularisants, et
\begin{equation}
h_\varepsilon= f \ast K_\varepsilon
\end{equation}
au sens
\begin{equation*}
h_\varepsilon(e^{it}) = \int f(e^{i(t-s)}) K_\varepsilon(e^{is})ds
\end{equation*}
en écrivant, ici comme dans la suite, $ \int $ à la place de $ \frac{1}{2\pi} \int_{-\pi}^\pi$. Posons
\begin{equation}
h_\varepsilon =\rho_\varepsilon e^{i\varphi_\varepsilon},\ \ 0\le \rho_\varepsilon <1\,.
\end{equation}
Quand $ \varepsilon\rightarrow 0 $, les normes $ \parallel f- h_\varepsilon\parallel_\infty $, $\parallel 1-\rho_\varepsilon\parallel_\infty$ et $ \parallel e^{i\varphi_\varepsilon}-f\parallel_\infty $ tendent vers 0. Donc, pour $ \varepsilon $ assez petit,
\begin{equation}
\deg e^{i\varphi_\varepsilon}=0\,.
\end{equation}
Considérons maintenant la fonction anti--analytique
\begin{equation}
R_\varepsilon =\exp(-\log \rho_\varepsilon + i \Hh(\log \rho_\varepsilon))\,,
\end{equation}
où $ \Hh $ est la transformation de Hilbert, et posons
\begin{equation}
H_\varepsilon = h_\varepsilon R_\varepsilon = \exp (i \varphi_\varepsilon + i \Hh (\log \rho_\varepsilon))\,.
\end{equation}
Comme les $ K_\varepsilon $ sont régularisants, $H_\varepsilon \in C^\infty (S^1,S^1)$.

La transformation de Hilbert applique $ L^\infty $ dans $ BMO $, et on sait que dans $ VMO(S^1,S^1) $ le degré est une fonction continue en norme $ BMO$([3], ou [1], theorem~1). D'où
\begin{equation}
\begin{matrix}
\parallel \exp(i \Hh(\log \rho_\varepsilon))\parallel_{BMO} \le \parallel\Hh (\log \rho_\varepsilon)\parallel_{BMO}\\
\noalign{\vskip2mm}
 \le \parallel\log \rho_\varepsilon\parallel_\infty =o(1) \quad (\varepsilon \rightarrow 0)
\end{matrix}
\end{equation}
 donc
\begin{equation}
\deg \exp (i\, \Hh(\log \rho_\varepsilon)) = 0
\end{equation}
quand $ \varepsilon $ est assez petit. Compte tenu de (13) et (15), cela donne
\begin{equation}
\deg H_\varepsilon =0\,.
\end{equation}
On peut alors appliquer la formule de Brézis (1), avec $ H_\varepsilon $ au lieu de $ f $, ce qui donne
\begin{equation}
\sum_{-\infty}^\infty |n|\, |\hat{H}_\varepsilon(n)|^2 = 2 \sum_0^\infty n |\hat{H}_\varepsilon(n)|^2\,.
\end{equation}
Estimons le second membre de (19). Pour $ I\subset \Z $, posons
\begin{equation}
(P_IF)^\wedge = \hat{F} \1_I\,.
\end{equation}
Comme $ R_\varepsilon $ est anti--analytique, on a d'après (15)
\begin{equation}
P_{[2^{k-1},2^k]} (H_\varepsilon) = P_{[2^{k-1},2^k]} (R_\varepsilon P_{[2^{k-1},\infty[} (h_\varepsilon))\,.
\end{equation}
Prenons les normes dans $ L^2  $ :
\begin{equation}
\begin{matrix}
\|P_{[2^{k-1},2^k]} (H_\varepsilon)\parallel_2 &\le \|R_\varepsilon P_{[2^{k-1},\infty[} (h_\varepsilon)\|_2\hfill\\
\noalign{\vskip2mm}
&\le \| \dfrac{1}{\rho_\varepsilon}\|_\infty \|P_{[2^{k-1},\infty[} (h_\varepsilon)\|_2\\
\noalign{\vskip2mm}
&\le 2 \|P_{[2^{k-1},\infty[} (f)\|_2\hfill\\
\noalign{\vskip2mm}
&= 2\Big(\displaystyle\sum_{n\ge 2^{k-1}}|a_n|^2\Big)^{1/2} \hfill
\end{matrix}
\end{equation}
(on a supposé, ce qui est permis, $ \int K_\varepsilon=1 $ ou proche de 1). Comme
\[ 
\sum_0^\infty n |\hat{H}_\varepsilon (n)|^2 \le \sum_{k\ge 0} 2^k \|P_{[2^{k-1},2^k]} (H_\varepsilon)\|_2^2\,,
 \]
 (22) donne
 \begin{equation}
 \begin{array}{ll}
\displaystyle\sum_0^\infty n|\hat{H}_\varepsilon(n)|^2 &\le \displaystyle\sum_{k\ge 0} 2^{k+2} \sum_{n\ge 2^{k-1}} |a_n|^2\\
\noalign{\vskip2mm}
&\le 16\ \displaystyle\sum_{n\ge 0} n|a_n|^2\\
\noalign{\vskip2mm}
&\le 16\ C
\end{array}
\end{equation}
d'après l'hypothèse (7). En vertu de (19), nous avons la borne uniforme en $ \varepsilon $ (pour $ \varepsilon<\varepsilon_0 $ assez petit pour avoir (17) donc (18))
\begin{equation}
\sum_{-\infty}^\infty |n| \ |\hat{H}_\varepsilon(n)|^2 \le 32\ C\,.
\end{equation}
Comme 
\[ 
\|f-H_\varepsilon\|_2 \le \|f-e^{i\varphi_\varepsilon}\|_2 + \|\Hh (\log \rho_\varepsilon)\|_2
 \]
 d'après la définition de $ H_\varepsilon $ en (15), et que le second membre tend vers 0 quand $ \varepsilon\rightarrow 0 $, on obtient à partir de (24)
 \begin{equation}
\sum_{-\infty}^\infty |n|\ |a_n|^2 \le 32\ C\,,
\end{equation}
la conclusion souhaitée en (8), avec $ C'=32\ C $.

\vspace{2mm}

\textbf{Cas} $ s >\dfrac{1}{2} $

Soit $ g_s(x)  $ $ (x\in\R )$ une fonction continue paire portée par l'intervalle $ [-2,2] $, positive et de classe $ C^2 $ sur l'intervalle ouvert $ ]-2,2[ $, et égale à $ (2-|x|)^{2s} $ quand $ 1\le |x| \le 2 $. Posons
\begin{equation}
K_{N,s}(t) = \sum_{n} g_s\big(\dfrac{n}{N}\Big) e^{int}\,.
\end{equation}

\vspace{2mm}

\textsc{Lemme}.-- \textit{On a}
\begin{equation}
|K_{N,s}(t)| \le c N(1 \wedge (N\|t\|)^{-1-2s})
\end{equation}
avec $ c=c_s $ et $ \|t\| = \mathrm{dist} (t,2\pi \Z) $.

\vspace{2mm}

\textbf{Preuve.} On vérifie en intégrant par parties deux fois que
\[ 
\int_\R e^{ixu}g_s(x)dx \le c (1\wedge |u|^{-1-2s})
 \]
 où $ c $ ne dépend que de $ s $ et de la ligne écrite, d'où
 \[ 
 \int_\R e^{ixu}g_s\big(\dfrac{x}{N}\big)dx \le c N(1\wedge |Nu|^{-1-2s})
  \]
  et la formule de Poisson donne (27).
  
  On pose
  \begin{equation}
\Delta_{N,s}(n) = N^{2s} g\big(\dfrac{n}{N} -2\big)\,.
\end{equation}
Ainsi
\begin{equation}
\begin{array}{ll}
\Delta_{N,s}(n) = n^{2s} &\mathrm{quand }\ 0\le n \le N\\
\noalign{\vskip2mm}
\Delta_{N,s}(n) \ge 0 &\mathrm{pour\ tout }\ n\,.
\end{array}
\end{equation}
\begin{equation}
\sum \Delta_{N,s}(n) e^{int} = N^{2s} e^{2iNt}K_{N,s}(t)\,.
\end{equation}
Nous allons calculer et comparer
\begin{equation}
\left \{
\begin{array}{l}
I_N=\sum(\Delta_{N,s}(n) +\Delta_{N,s}(-n))|a_n|^2\\
\noalign{\vskip2mm}
J_N=\sum (\Delta_{N,s}(n)-\Delta_{N,s}(-n))|a_n|^2\,.
\end{array}\right .
\end{equation}

De nouveau nous supposons (2), (7) et $ \deg f=0 $, et nous nous proposons d'établir (8). Comme $ \deg f=0 $, nous pouvons écrire
\begin{equation}
f(e^{it}) = e^{i\varphi(t)}\,,\ \ \varphi\in C(\R/2\pi\Z, \R)\, .
\end{equation}
\begin{equation}
a_n = \int e^{i\varphi(t)}e^{-int}dt
\end{equation}
et
\begin{equation}
|a_n|^2 = \intint e^{i(\varphi(t_1)-\varphi(t_2))}e^{-in(t_1-t_2)}dt_1 dt_2
\end{equation}
puis, suivant (30)
\begin{equation}
\begin{array}{l}
\displaystyle\sum \Delta_{n,s}(n)|a_n|^2 \!=\!\! \intint e^{i(\varphi(t_1)-\varphi(t_2))} N^{2s} e^{-2iN(t_1-t_2)} K_{N,s}(t_1-t_2) dt_1dt_2\\
\noalign{\vskip2mm}
\displaystyle\sum \Delta_{n,s}(-n)|a_n|^2\! = \!\!\intint e^{i(\varphi(t_1)-\varphi(t_2))} N^{2s} e^{-2iN(t_1-t_2)} K_{N,s}(t_1-t_2) dt_1dt_2
\end{array}
\end{equation}
et, suivant (31), en observant que $ I_N $ et $ J_N $ sont réelles,
\begin{equation}
\kern-0,5cm\left \{
\begin{array}{l}
\!\!I_N=2\!\displaystyle\intint \cos (\varphi(t_1)\!-\!\varphi(t_2)) N^{2s} \cos\ 2N (t_1\!-\!t_2)K_{N,s}(t_1\!-\!t_2)dt_1dt_2\\
\noalign{\vskip2mm}
\!\!J_N=2\!\displaystyle\intint \sin (\varphi(t_1)\!-\!\varphi(t_2)) N^{2s} \sin\ 2N (t_1\!-\!t_2)K_{N,s}(t_1\!-\!t_2) dt_1dt_2\,.
\end{array}\right .
\end{equation}
Désormais seul $ J_N $ nous importe. Comme la valeur moyenne de $ \sin\, 2Nt\, K_{N,s}(t) $ est nulle, on peut écrire
\begin{equation}
\begin{array}{c}
J_N=2\displaystyle\intint( \sin (\varphi(t_1)-\varphi(t_2))-(\varphi(t_1)-\varphi(t_2))) N^{2s}\\ \sin\ 2N (t_1-t_2)K_{N,s}(t_1-t_2)dt_1dt_2
\end{array}
\end{equation}
d'où
\begin{equation}
|J_N| \le \intint |\varphi(t_1)-\varphi(t_2)|^3 N^{2s} K_{N,s} (t_1-t_2)dt_1dt_2
\end{equation}
et, tenant compte de (27),
\begin{equation}
|J_N| \le c \intint |\varphi(t_1)-\varphi(t_2)|^3 (N^{1+2s}\wedge \|t_1-t_2\|^{-1-2s}) dt_1dt_2\,.
\end{equation}

L'hypothèse (7), compte tenu de (29), donne
\begin{equation}
\sum \Delta_{N,s}(n) |a_n|^2 \le C
\end{equation}
soit, dans la notation (31), $ I_N+J_N \le 2C $, et
\begin{equation}
\sum_{|n|\le N} |n|^{2s}|a_n|^2 \le I_N\le |J_N|+2C\,.
\end{equation}
Il s'agit de montrer que $\sup\limits_N |J_N|$ est fini.

Pour utiliser (39) et (41), nous allons, enfin, utiliser l'hypothèse $ s>\frac{1}{2} $. Jointe à (4), cela entra\^{i}ne $ \sum |n|\,|a_n|^2 <\infty $, soit
\begin{equation}
\sum_k 2^k s_k^2 <\infty \,, \quad s_k^2 = \sum_{2^k \le |n|<2^{k+1}}|a_n|^2\,.
\end{equation}
Prenons $ N=2^{j+1} $ tel que
\begin{equation}
2^j s_j^2 = \max_{k\ge j} 2^k s_k^2\,.
\end{equation}
Alors
\begin{equation}
\sum_{|n|\ge N}|a_n|^2 = \sum_{k> j} s_k^2 \le  s_j^2 =  \sum_{\frac{N}{2}\le |n|<N} |a_n|^2
\end{equation}
\begin{equation}
\begin{array}{c}
\displaystyle\sum_{n\in \Z}\big(|n|\wedge N\big)^{2s} |a_n|^2 \le \sum_{|n|<N} |n|^{2s} |a_n|^2 + N^{2s} \sum_{\frac{N}{2}\le |n|<N}|a_n|^2\\
\noalign{\vskip2mm}
\displaystyle \le (1+2^{2s}) \sum_{|n|\le N} |n|^{2s}|a_n|^2\,.
\end{array}
\end{equation}

Dans ce qui suit écrivons $c'\ (\mathrm{resp}\ C')$ pour un nombre $ >0 $ qui pourra dépendre de $s \ (\mathrm{resp}\ C,s,f) $ et de la ligne écrite, mais non de $ N $. Ainsi, on a
\begin{equation}
\begin{array}{c}
\displaystyle\sum_{n\in \Z}(|n|\wedge N)^{2s} |a_n|^2\ge c' \intint |f(e^{it_1})- f(e^{it_2})|^2 \\
\noalign{\vskip2mm}
(N^{1+2s} \wedge \|t_1-t_2\|^{-1-2s}) dt_1dt_2
\end{array}
\end{equation}
où $c'$ est une constante absolue comme le montre un calcul du second membre en fonction des $a_n$ (c'est ici qu'intervient l'hypothèse $s<1$) et, d'après (39), (41) et (45), 
\begin{equation}
\begin{array}{c}
\displaystyle\intint |f(e^{it_1})-f(e^{it_2})|^2 (N^{1+2s} \wedge \|t_1-t_2\|^{-1-2s}) dt_1dt_2\\
\noalign{\vskip2mm}
\displaystyle\le C'+c' \intint |\varphi(t_1)-\varphi(t_2)|^3 (N^{1+2s}\wedge \|t_1-t_2\|^{-1-2s})dt_1dt_2\,.
\end{array}
\end{equation}

Comme $ \varphi\in C(\R/2\pi\Z,\, \R)$, il existe pour tout $ \beta >0$
 donné un $ \alpha=\alpha(\varphi,\beta)>0 $ tel que $ \|t_1-t_2\|<\alpha $ entra\^{i}ne $ \varphi(t_1)-\varphi(t_2)<\beta $, et
 \[ 
|f(e^{it_1}) -f(e^{it_2})| = |e^{i\varphi(t_1)}-e^{i\varphi(t_2)}| \simeq |\varphi(t_1)-\varphi(t_2)|\,.
  \]
  Alors
  \begin{equation}
\begin{array}{c}
\displaystyle\intint |\varphi(t_1)-\varphi(t_2)|^3 (N^{1+2s}\wedge \|t_1-t_2\|^{-1-2s})dt_1dt_2\\
\noalign{\vskip2mm}
\displaystyle = \intint_{\|t_1-t_2\|\ge\alpha} +\intint_{\|t_1-t_2\|<\alpha}\\
\noalign{\vskip2mm}
\displaystyle \le c' \alpha^{-2s} \|\varphi\|_\infty^3 +c' \beta \!\intint\! |f(e^{it_1})\\
\noalign{\vskip2mm}
\hfill-f(e^{it_2})|^2 (N^{1+2s}\wedge \|t_1-t_2  \|^{-1-2s})dt_1dt_2\,.
\end{array}
\end{equation}
En choisissant $ \beta $ assez petit, on obtient d'après (47) et (48)
\begin{equation}
\intint |f(e^{it_1}) -f(e^{it_2})|^2 (N^{1+2s} \wedge \|t_1-t_2  \| ^{-1-2s})dt_1dt_2 <C'
\end{equation}
et $ C' $ ne dépend pas de $ N $, d'où résulte $ f\in H^s(S^1,S^1) $, c'est--\`{a}--dire (8).

\vskip2mm

\textbf{Cas} $ s<\frac{1}{2} $.

Ici tout ce qui suit (41) est en défaut. Mais tout ce qui va de (26) à (41) est valable, à commencer par le lemme et la formule (27), qui s'établit ici avec une seule intégration par parties. Nous allons utiliser ces formules en remplaçant $ f $ par la fonction $ H_\varepsilon $ définie en (15).

Auparavant, quitte à changer $ C $ dans (7), nous nous ramenons au cas où l'argument de $ f $ est proche de~0,
\begin{equation}
|\arg f| < \delta_0
\end{equation}
($ \delta_0 $ sera fixé plus tard et dépendra seulement de la constante $ c=c_0 $ figurant en (39)) ; il suffit pour cela de multiplier la fonction $ f $ donnée par une fonction $ \in C^{\infty} (S^1,S^1)$ convenable. Il en résulte, dans les notations (11) et (12), que
\begin{equation}
|\varphi_\varepsilon| < \delta_0
\end{equation}
(en se restreignant aux noyaux de convolution $ K_\varepsilon $ de moyenne 1). Ecrivons, suivant (15),
\begin{equation}
\begin{array}{l}
H_\varepsilon=e^{i\varphi}\\
\varphi=\varphi_\varepsilon + \Hh(\log \rho_\varepsilon)
\end{array}
\end{equation}
(c'est la nouvelle signification de $ \varphi $). D'après (16) et (51), on peut choisir $ \varepsilon $ assez petit pour que
\begin{equation}
\| \varphi \|_{BMO}<\delta_0\,.
\end{equation}

De nouveau nous avons $ H_\varepsilon \in C^\infty (S^1,S^1) $. Nous allons d'abord établir l'analogue de (7) pour $ H_\varepsilon $, à savoir
\begin{equation}
\sum_0^\infty n^{2s} |\hat{H}_\varepsilon(n)|^2 <C\,,
\end{equation}
$ C $ désignant ici et à partir de maintenant un nombre $ C(f,s) $ dépendant de $ f $ et de $ s $, mais non de $ \varepsilon $. Il suffit de reprendre les calculs de (20) à (23), en écrivant ici
\begin{equation}
\begin{array}{ll}
\displaystyle\sum_0^\infty n^{2s}  |\hat{H}_\varepsilon(n)|^2 &\le \displaystyle\sum_{k\ge 0} 4^{ks} \| P_{[2^{k-1},2^k]} (H_\varepsilon) \|_2^2\\
\noalign{\vskip2mm}
&\displaystyle \le\sum_{k\le 0} 4^{ks+1} \sum_{n\ge 2^{k-1}} |a_n|^2\\
\noalign{\vskip2mm}
&\displaystyle \le 4^{s+1} (1-4^{-s})^{-1} \sum_{n\ge 0} |a_n|^2
\end{array}
\end{equation}
qui est bien de la forme (54).

Partant de (52), suivons les notations et les calculs de (26) à (41), en remplaçant $ a_n $ par $ \hat{H}_\varepsilon(n) $. Comme $ \varphi\in C^\infty $, (39) et (41) donnent
\begin{equation}
|J|\le c_0 \intint \frac{|\varphi(t_1)-\varphi(t_2)|^3}{\| t_1-t_2  \|^{1+2s}} dt_1 dt_2
\end{equation}
puis, en faisant tendre $ N $ vers l'infini
\begin{equation}
\sum_{n\in \Z} |n|^{2s} |\hat{H}_\varepsilon(n)|^2 \le c_0 \intint \frac{|\varphi(t_1)-\varphi(t_2)|^3}{\|t_1-t_2 \|^{1+2s}} dt_1dt_2 +2C
\end{equation}
$ C $ étant le second membre de (54). Le premier membre de (57) est $ \| H_\varepsilon \| ^2_{H^s}$, qui, à une équivalence numérique près, s'écrit aussi sous forme d'intégrale. L'inégalité
\begin{equation}
|u-v| \le |e^{iu}-e^{iv}| + |u-v|^{3/2} \quad (u,v\in \R)
\end{equation}
permet d'écrire
$$
\begin{array}{c}
\displaystyle\intint \frac{|\varphi(t_1)-\varphi(t_2)|^2}{\|t_1-t_2 \|^{1+2s}} dt_1dt_2 \\
\noalign{\vskip2mm}
\displaystyle\le \intint \frac{|H_\varepsilon(e^{it_1})- H_\varepsilon(e^{it_2})|^2}{\|t_1-t_2 \|^{1+2s}} dt_1dt_2 + 
\intint \frac{|\varphi(t_1)-\varphi(t_2)|^3}{\|t_1-t_2 \|^{1+2s}} dt_1dt_2
\end{array}
$$
d'où, par (57),
\begin{equation}
\intint \frac{|\varphi(t_1)-\varphi(t_2)|^2}{\|t_1-t_2 \|^{1+2s}} dt_1dt_2 \le (Ac_0+1) \intint \frac{|\varphi(t_1)-\varphi(t_2)|^3}{\|t_1-t_2 \|^{1+2s}} dt_1dt_2 +AC
\end{equation}
$ A$  étant une constante absolue. Les deux membres de (59) s'expriment à l'aide des $ P_{[2^{k-1},2^k]}(\varphi) $ (définis comme en (20)), et (59) s'écrit
\begin{equation}
\sum_{k\ge 0} 4^{ks} \|P_{[2^{k-1},2^k]}(\varphi)\|_2^2 \le B(c_0+1) 
 \sum_{k> 0} 4^{ks} \|P_{[2^{k-1},2^k]}(\varphi)\|_3^3 +  BC
 \end{equation}
$ B $ étant une constante absolue. Or
\begin{equation}
\|P_{[2^{k-1},2^k]}(\varphi)\|_\infty \le \| \varphi\|_{BMO}
\end{equation}
donc, d'après (53)
\begin{equation}
\|P_{[2^{k-1},2^k]}(\varphi)\|_3^3 \le \delta_0 \|P_{[2^{k-1},2^k]}(\varphi)\|_2^2\,.
\end{equation}
Choisissons au départ $ \delta_0= \frac{1}{2B(c_0+1)} $. Alors, d'après (60) et (62),
\begin{equation}
\|\varphi\|_{H^s}^2 \le C=C(f,s)
\end{equation}
et il en résulte
\begin{equation}
\|H_\varepsilon\|_{H^s}^2 \le C=C(f,s)
\end{equation}
uniformément par rapport à $ \varepsilon $, et en faisant tendre $ \varepsilon $ vers 0 on a 
\begin{equation}
f\in H^s(S^1,S^1)\,,
\end{equation}
la conclusion voulue.

\vskip2mm

\textbf{Récapitulation et preuve du théorème 2}.

La preuve du théorème 1 s'est déroulée en trois étapes, et les principaux ingrédients sont apparus dans l'examen des cas $ s=\frac{1}{2} $ et $ s>\frac{1}{2} $. Ces ingrédients sont utilisés dans le cas $ s<\frac{1}{2} $, dont le traitement s'étend immédiatement au cas général $ 0 < s < 1 $. Pour établir le théorème 2, il suffit donc d'adapter la preuve donnée dans ce dernier cas en remplaçant l'hypothèse $ f\in C(S^1,S^1) $ par $ f\in VMO(S^1,S^1) $.

Seul le début est à changer, jusqu'à la formule (53). Tout ce qui suit (53) est à conserver littéralement. 

A la place de (50), nous nous ramenons au cas $ \deg f=0 $ et
\begin{equation}
\| f  \|_{BMO} < \delta_0\,,\ \| \arg f \|_{BMO} <\delta_0
\end{equation}
d'où résulte
\begin{equation}
\| \varphi _\varepsilon \|_{BMO} < \delta_0\,.
\end{equation}
D'après Brézis et Nirenberg [4] on a
\begin{equation}
\|1-\rho_\varepsilon\|_\infty =o(1) \quad (\varepsilon\rightarrow 0)\, ;
\end{equation}
rappelons la preuve :
\begin{equation}
1 - \Big| \frac{1}{|I|} \int_I f\Big| \le \frac{1}{|I|} \int_I \Big| f-\frac{1}{|I|} \int f\Big| =o(1) \quad (|I|\rightarrow0)
\end{equation}
d'après (66), d'où $ 1-|h_\varepsilon| =o(1) $ $ (\varepsilon\rightarrow0) $, c'est--à--dire (68). (67) et (68) donnent
\begin{equation}
\| \varphi_\varepsilon+\Hh(\log \rho_\varepsilon) \|_{BMO} < \delta_0+o(1) \quad (\varepsilon\rightarrow 0)\,,
\end{equation}
ce qui établit (53), d'où la conclusion.

\vskip4mm

Les auteurs remercient le referee par sa lecture attentive de
l'article, et en particulier pour une amélioration qu'il a suggérée, à
savoir de ne pas nous limiter au cas  $0 < s < 1$, mais de traiter le cas $s > 0$. 

\vskip4mm

\textbf{Réponse à la question Q3, extension au cas $s\ge 1$.}

\vskip2mm

\textsc{Théorème} 3.-- \textit{Les théorèmes $1$ et $2$ sont valables en remplaçant $0<s<1$ pour $s>0$.}

\vskip2mm

\textit{Preuve}. Soit $s\ge 1$, $f\in VMO(S^1,S^1)$ et $Pf\in H^s(S^1,\C)$, où $f\sim \sum\limits_{-\infty}^\infty a_nz^n$ et $Pf = \sum\limits_0^\infty a_nz^n$. Nous savons déjà que $f\in H^{s'}(S^1,S^1)$ pour tout $s'<1$, donc $f=e^{i\varphi}$ où $\varphi \in H^{s'}(S^1,\R)$ [1]. Quitte à multiplier $f$ par une fonction
de classe  $C^\infty$ comme nous l'avons fait dans l'étude du cas   $s  <  1/2$, on peut
 supposer $\|\varphi\|_{H^{s'}} < \delta<\frac{1}{10}$. Ecrivons $\|\varphi\| = \Big(\sum\limits_{-\infty}^\infty |\widehat{\varphi}_n|^2 |n|^{2s'}\Big)^{1/2}$. Alors 
\begin{eqnarray*}
f &=&1+ i\varphi +h,\ \|h\| \le \dfrac{1}{2} \|\varphi\|^2 + \dfrac{1}{6} \|\varphi\|^3 + \cdots \le \delta\|\varphi\|\\
\noalign{\vskip2mm}
Pf&=&1+i\, P\varphi+Ph,\ \|Ph\| \le \|h\| \le \delta\|\varphi\|
\end{eqnarray*}
et, comme $\varphi$ est réelle, $\|P\varphi\|^2 = \frac{1}{2}\|\varphi\|^2$. Donc
\begin{eqnarray}
\|f\| &\le& \|\varphi\| (1+\delta)\,, \nonumber\\ 
\noalign{\vskip2mm}
\|Pf\| &\ge& \dfrac{1}{2}\|\varphi\| -\delta\|\varphi\| \,, \nonumber\\
\noalign{\vskip2mm}
\|f\| &\le& \|Pf\| \dfrac{1+\delta}{\frac{1}{2}-\delta} \le 3\|Pf\|\,.
\end{eqnarray}
Sous la forme
$$
\sum_{-\infty}^\infty |n|^{2s'} |a_n|^2 \le 3 \sum_0^\infty \, n^{2s'} |a_n|^2\,,
$$
il est clair que (71) s'étend à toutes les valeurs de $s'$ pour lesquelles le second membre est fini. Comme c'est le cas pour $s'=s$ par hypothèse, on a bien $f\in H^s (S^1,S^1)$, et, comme cela est conservé par multiplication par une fonction appartenant à $C^\infty (S^1,S^1)$, le théorème est établi.

\eject

\vskip4mm

\parindent=0mm

Jean Bourgain

bourgain@math.ias.edu

Institute of Advanced Study,

Princeton,

NJ, USA

\vskip4mm
J.--P.  Kahane

jean-pierre.kahane@math.u-psud.fr

Laboratoire de Mathématiques,

Université Paris--Sud,

91405 Orsay Cedex

France

\end{document}